\documentclass{article}

\usepackage{amssymb}
\usepackage{amsmath}
\usepackage{amsthm}
\usepackage{a4wide}

\usepackage[english]{babel}

\newcommand{\clin}{\centerline}

\newcommand{\ovline}{\overline}
\newcommand{\wtilde}{\widetilde}

\newcommand{\vect}[1]{\mathbf{#1}}

\newcommand{\larr}{\leftarrow}
\newcommand{\rarr}{\rightarrow}

\newcommand{\Rarr}{\Rightarrow}
\newcommand{\LRarr}{\Leftrightarrow}

\newcommand{\lesq}{\leqslant}
\newcommand{\greq}{\geqslant}

\newcommand{\veps}{\varepsilon}
\newcommand{\vphi}{\varphi}

\newcommand{\forceparindent}{\hskip 1.5em}

\DeclareMathOperator{\sign}{sign}

\newcounter{Formulanum}
\newcommand{\equ}[1]{\refstepcounter{Formulanum}%
\label{#1}%
(\arabic{Formulanum})}
\newcommand{\equref}[1]{(\ref{#1})}

\newcounter{Propositionum}
\newcommand{\prop}[1]{\refstepcounter{Propositionum}%
\label{#1}%
\textbf{Proposition \arabic{Propositionum}.}}
\newcommand{\propref}[1]{Prop. \ref{#1}}

\newcounter{Examplnum}
\newcommand{\examp}[1]{\refstepcounter{Examplnum}%
\label{#1}%
\textbf{Example \arabic{Examplnum}.}}
\newcommand{\exampref}[1]{Ex. \ref{#1}}

\newenvironment{mproof}{\textbf{Proof.}}{\qed}

\newenvironment{mexamp}{$\vartriangleleft$}{$\vartriangleright$}

\begin{document}

\title{On sound ranging in Hilbert space}

\author{Sergij V. Goncharov\thanks{Faculty of Mechanics and Mathematics, Oles Honchar Dnipro National University, 72 Gagarin Avenue, 49010 Dnipro, Ukraine.
\textit{E-mail: goncharov@mmf.dnulive.dp.ua}}}

\date{July 2017}

\maketitle

\begin{abstract}
We consider the sound ranging problem, which is to find the position of the source-point from the moments when the wave-sphere of linearly, with time, increasing radius
reaches the sensor-points, in the infinite-dimensional separable Euclidean space H, and describe the solving methods,
for entire space and for its unit sphere. In the former case, we give some sufficient conditions for uniqueness of the solution.
We also provide two examples with the sets of sensors being a basis of H:
1st, when sound ranging problem and so-called dual problem both have single solutions,
and 2nd, when sound ranging problem has two distinct solutions.
\end{abstract}

\clin{\small \textit{MSC2010:} Pri 46C05, Sec 40A05, 15A15}

\clin{\small \textit{Keywords:} Hilbert space, infinite dimensional, sound ranging, TDOA, sphere, basis, uniqueness, series}

\newlength{\oldlineskip}
\setlength{\oldlineskip}{\baselineskip}
\baselineskip=9pt
\makeatletter
\@starttoc{toc} 
\makeatother
\setlength{\baselineskip}{\oldlineskip}

\section*{Introduction}
\addcontentsline{toc}{section}{Introduction}

\forceparindent
By \textit{sound ranging} (SR), we mean the following problem. Let $(X;\rho)$ be a metric space, i.e. the set $X$ with metric $\rho \colon X\times X \rarr \mathbb{R}_+$.
Let $\vect{s} \in X$ be an unknown point, ``source''. At unknown moment $t_e \in \mathbb{R}$ of time the source ``emits the (sound) wave'',
which is the sphere

\clin{$S\bigl(\vect{s}; v(t-t_e)\bigr) = \bigl\{ \vect{x} \in X \mid \rho(\vect{x}; \vect{s}) = v(t - t_e) \bigr\}$}

\noindent
for any moment $t \greq t_e$. Here $v$ is known ``sound velocity'', and we may assume, without loss of generality, that $v = 1$
(switching to scaled time $t \larr vt$ if $v \ne 1$).

Let $R = \{ \vect{r}^{(i)} \}_{i \in I}$, $\vect{r}^{(i)} \in X$, be an indexed set of ``sensors'', whose positions are known.
Suppose that for each sensor we know the moment $t_i$ when it was reached by the expanding wave;
that is, $t_i = t_e + \rho (\vect{r}^{(i)};\vect{s})$ are known.

The problem is to find $\vect{s}$ and $t_e$, --- from known moments when the wave reaches known sensors, $(R;\{ t_i \})$.
We're also interested in uniqueness of the solution.

{\ }

\textbf{Retrospection.} The researches on SR, ---
also called \textit{passive location}, \textit{sound triangulation}, \textit{time-(difference-)of-arrival source localization}, ---
are considered to begin at the times of World War I, with the works of William L. Bragg, Lucien Bull, Erich Waetzmann among the others
(see \cite{bazzoni1918}, \cite{bragg1921}, \cite{bull1915}, and \cite{macleod2000}, \cite{schiavon2015} for a survey).
Similar questions were studied in \cite{bateman1918}, \cite{dadourian1921}, \cite{esclangon1925}, \cite{hopejones1928}
from that ``geometrical'' era. During the century that followed, along with military
(\cite{canistraro1996}, \cite{ferguson2002}, \cite{hercz1987}, \cite{miller1976}, \cite{rachele1966})
and surveillance-related (\cite{maroti2004}) applications,
SR attracted acousticians and seismologists, to name a few, --- \cite{friedlander1987}, \cite{liu2016}, \cite{schriever1952}, \cite{watkins1971}, \cite[5.7]{whitlow2008}.
For a fairly long time, SR problems accompany the studies of (wireless) sensor networks
(\cite{bancroft1985}, \cite{burgess2013}, \cite{duraiswami1999}, \cite{lombard2006}, \cite{mao2009}, \cite{schloemann2014}, \cite{shen2012}; see also \cite{zekavat2012}).

{\ }

Naturally, the majority of these researches relates to $\mathbb{R}^2$ or $\mathbb{R}^3$
(emphasized occasionally: \cite{tobias1976}, \cite{watkins1971}, \cite{wu2010}),
though there are exceptions (\cite{burgess2013}, \cite{lombard2006}, \cite{robinson2012}, \cite{schloemann2014}, \cite{torrieri1984}).
The basic problem was generalized to take into account the ``wind'' or ``flow'' that gives an additional movement to the wave,
the variations of sound velocity at different regions of space, diffraction and reverberation,
the inaccuracy of measurements and its influence on the solution(s),
noise removal etc. (\cite{hoock1979}, \cite{lee1973}, \cite{liu2016}), --- generally speaking, the factors imposed by ``physics''
(as a result, sometimes there's the inclination towards practical applicability instead of rigour).

The uniqueness of the solution was analyzed as well, for example in \cite{compagnoni2014}, \cite{liu2016}, \cite{robinson2012}, \cite{spencer2010}, \cite{spiesberger2001}.

The substantial part of the studies in the field deals with the so-called overdetermined problems, when the data from each sensor contains some random error,
and the position of the source is estimated in attempt to reduce the uncertainty
(\cite{burgess2013}, \cite{duraiswami1999}, \cite{friedlander1987}, \cite{gillette2008}, \cite{lee1973}, \cite{lee1967}, \cite{torrieri1984}).

{\ }

\textbf{Aim.} The generalization being investigated here concerns only the infinite dimensionality of the (empty) space where source and sensors are placed,
and omits ``physical'' factors (comparing it with most of the papers listed in References, we warn at once that it is of much less ``applicability'').

{\ }

Consider (associated, in a way) question about solving such a problem: what limitations do we impose on the ``procedure'' --- or ``algorithm'' ---
of obtaining the solution $\vect{s} \in X$?
There's always a ``universal'' one, $\mathcal{U}$: ``go over all $\vect{s} \in X$ and select those satisfying $t_i = t_e + \rho (\vect{r}^{(i)};\vect{s})$''
(if we know $\vect{s}$, $t_e$ can be found easily as $t_i - \rho (\vect{r}^{(i)}; \vect{s} )$ for some $i$).
But $\mathcal{U}$ seems to be too ``heavy''. If, in a sense, the verification of each $\vect{s}$ takes a non-zero amount of time $\tau$,
then for many kinds of spaces $X$ we'll need an infinite, non-countable set of ``verificating entities'' to confine into a finite time.

Less heavy but still not appropriate procedure $\mathcal{N}$ is as follows: take some point $\vect{b} \in X$ as the ``origin'' and, at each $n$-th step,
make the ``$\frac{1}{n}$-net $N_n$'' in the ball $B(\vect{b};n)$ ($\forall \vect{x} \in B(\vect{b};n)$ $\exists \vect{y} \in N_n$: $\rho (\vect{x}; \vect{y}) < \frac{1}{n}$).
In turn, for each $\vect{y} \in N_n$ calculate the ``defect'' $\delta (\vect{y}) = \sup\limits_i \bigl| t_e + \rho (\vect{r}^{(i)};\vect{y}) - t_i \bigr|$,
and select the $\vect{y}$ with defect not greater than $\inf\limits_{\vect{y}\in N_n} \delta (\vect{y})+\frac{1}{n}$. $\rho (\vect{b}; \vect{s}) < \infty$ implies
the existence of the sequence $\{ \vect{y}_n' \in N_n \}$ such that $\vect{y}_n' \rarr \vect{s}$ as $n \rarr \infty$ (not unique, though).

We prefer more ``countable'' and ``aimed'' methods. On the other hand, we ``allow ourselves'' the calculation of infinite sums in a finite time, as short as we want.

{\ }

\textbf{Well-knowns.} Hereinafter, $H$ is \textsl{the} separable infinite-dimensional Hilbert space over the field of reals $\mathbb{R}$. We denote by ${<}\vect{x};\vect{y}{>}$
the scalar product of any $\vect{x}, \vect{y}\in H$, and by $\| \vect{x} \|$ the norm of $\vect{x}$. As usual, ${<}\vect{x};\vect{x}{>} = \| \vect{x} \|^2$.
Since the field is $\mathbb{R}$, ${<}\vect{y};\vect{x}{>} = {<}\vect{x};\vect{y}{>}$ (the complex case reduces to the real one with ``twice more dimensions'',
due to representability of distance between 2 points with complex coordinates through their real and imaginary parts).
$\vect{x} \perp \vect{y}$ means ${<}\vect{x};\vect{y}{>} = 0$.

Some common properties of scalar product and norm are used without explicit reference:

$\bullet$ for any orthonormal basis $\{ \vect{e_k} \}_{k\in \mathbb{N}}$ of $H$, if $\vect{x} = \sum\limits_{k = 1}^{\infty} x_k \vect{e_k}$
and $\vect{y} = \sum\limits_{k = 1}^{\infty} y_k \vect{e_k}$, then ${<}\vect{x};\vect{y}{>} = \sum\limits_{k = 1}^{\infty} x_k y_k$,
independent of basis $\{ \vect{e_k} \}$.

$\bullet$ Cauchy-Bunyakowsky-Schwartz inequality (CBS): $|{<}\vect{x};\vect{y}{>}| \lesq \| \vect{x} \| \cdot \| \vect{y} \|$,
and the inequality becomes equality if and only if $\vect{x}$ and $\vect{y}$ are linearly dependent, that is,
$\exists a, b \in \mathbb{R}$: $a^2 + b^2 \ne 0$ and $a \vect{x} + b \vect{y} = \theta$ (where $\theta$ is the zero of $H$ as linear vector space).
Moreover, if ${<}\vect{x};\vect{y}{>} = \| \vect{x} \| \cdot \| \vect{y} \|$ and $\vect{y} \ne \theta$, then $\vect{x} = c \vect{y}$, where $c \greq 0$.

$\bullet$ If $\| \vect{x} \pm \vect{y} \| = \| \vect{x} \| \pm \| \vect{y} \|$, then $\vect{x}$ and $\vect{y}$
are linearly dependent (follows from previous statement).

$\bullet$ $\bigl| \| \vect{x} \| - \| \vect{y} \| \bigr| \lesq \| \vect{x} - \vect{y} \|$.

{\ }

\textbf{Disclaimer.} It appears that $X=H$ introduces nuances (mostly dealing with limits and convergency),
however the basic method is that of $\mathbb{R}^n$ case.
\textsl{We surmise many of the subsequent results, --- especially the ``auxiliary'' ones, --- to be already known, even as ``folklore''}, perhaps;
if so, this is merely where they come together... once more (see also ``Acknowledgements'').

\section{SR in Hilbert space}

\forceparindent
Hereinafter, the set $R \subset H$ of sensors is finite or countable, and the SRP(roblem) is supposed to have at least one solution $(\vect{s}_0;t_{e;0})$, which may be unknown.

To simplify the notation, we move ``the origin of space and time'' to one of sensors at the moment when the (sound) wave reaches it.
So, $R = \{ \vect{r}^{(0)}, \vect{r}^{(1)}, ..., \vect{r}^{(n)}, ... \}$ with $\vect{r}^{(0)} = \theta$, and the wave reaches these sensors
at the moments $t_0 = 0$, $t_1$, $t_2$, ..., where $t_i = t_{e;0} + \| \vect{r}^{(i)} - \vect{s}_0 \|$.

By $L(A)$ we denote the linear closure of the set $A \subseteq H$. Note that $L(R) = L(\{ \vect{r}^{(i)} \}_{i\in \mathbb{N}})$, of all sensors but $\vect{r}^{(0)}$.
We denote $\{ \vect{r}^{(i)} \}_{i\in \mathbb{N}}$ by $\dot{R}$.

{\ }

We begin by excluding the sets of sensors such that the solution, if it exists, is obviously not unique.
If $L(\dot{R}) \ne H$ and the source $\vect{s}_0 \notin L(\dot{R})$, then by projection theorem $\vect{s}_0 = \vect{u} + \vect{h}$, where $\vect{u} \in L(\dot{R})$,
$\vect{h} \perp L(\dot{R})$ and $\vect{h} \ne \theta$. Then for each sensor the square of ``reaching time'',

\clin{$(t_i - t_{e;0})^2 = \| \vect{r}^{(i)} - \vect{s}_0 \|^2 = \| \vect{r}^{(i)} - \vect{u} - \vect{h} \|^2 = \| \vect{r}^{(i)} - \vect{u} \|^2 + \| \vect{h} \|^2$}

\noindent
(${<}\vect{r}^{(i)} - \vect{u}; \vect{h}{>} = 0$ since these elements are orthogonal) would be the same for $\vect{s'_0} = \vect{u} - \vect{h}$, ---
SRP has 2 solutions being non-distinguishable by the $\{ t_i \}$.

Moreover, if $H = L(\dot{R}) \oplus K$ and $\dim K \greq 2$, then $\exists \vect{w} \ne \theta$: $\vect{w} \perp L(\dot{R})$ and $\vect{w} \perp \vect{h}$.
Consider normalized $\wtilde{\vect{h}} = \frac{\vect{h}}{ \| \vect{h} \|}$, $\wtilde{\vect{w}} = \frac{\vect{w}}{ \| \vect{w} \|}$,
and let $\vect{s}(\vphi) = \vect{u} + \| \vect{h} \| (\cos \vphi \cdot \wtilde{\vect{h}} + \sin \vphi \cdot \wtilde{\vect{w}})$;
it is easy to see that $\vect{s}(\vphi)$ is a solution of SRP for any $\vphi \in [0;2\pi)$, --- we have an infinite, non-countable set of solutions.

{\ }

Let $L(\dot{R}) = H$, and let $\dot{R}$ be a linearly independent set. In other words, let $\dot{R}$ be a basis of $H$.

We introduce the orthonormal basis $B = \{ \vect{e}_i \}_{i\in \mathbb{N}}$, derived from $\dot{R}$ by Gram-Schmidt orthogonalization.
With $B$, $H$ is $l_2$ and $\vect{r}^{(i)} = \sum\limits_{j=1}^{\infty} r^{(i)}_j \vect{e_j} = (r^{(i)}_1; r^{(i)}_2; ...; r^{(i)}_i; 0; 0; ...)$, where $r^{(i)}_i \ne 0$.

To find $\vect{s}$ is to find its coordinates $(s_1; s_2; ...)$.

{\ }

The SRP is equivalent to the following set of equations: $t_i = t + \| \vect{r}^{(i)} - \vect{s} \|$, $i \in \mathbb{Z}_+$ \hfill \equ{emeqset}

Note that $0$-th equation is actually $0 = t + \| \theta - \vect{s} \|$ $\LRarr$ $t = - \| \vect{s} \|$.

{\ }

(For instance, take $H = L_2[a;b]$, and let $f \in L_2[a;b]$ be an unknown function. Suppose that for each $i\in \mathbb{N}$ we know
$t_i = t + \bigl( \int\limits_a^b |f(x) - x^{i-1}|^2 dx \bigr)^{\frac{1}{2}}$, where $t = - \bigl( \int\limits_a^b f^2(x) dx \bigr)^{\frac{1}{2}}$ is unknown too.)

{\ }

We now proceed to the implied set of equations $\| \vect{r}^{(i)} - \vect{s} \|^2 = (t_i - t)^2$, $i \in \mathbb{Z}_+$ \hfill \equ{impeqset}

\noindent
which may have additional solutions $(\vect{s}; t)$. To distinguish them from those of \equref{emeqset},
we verify that $\forall i\in \mathbb{Z}_+$: $t \lesq t_i$ (in particular, $t \lesq t_0 = 0$) --- ``the wave was emitted \textit{before} it reached sensors''.

{\ }

\textbf{Dual problem.} The additional solutions of \equref{impeqset} such that $t \greq t_i$ are the solutions of the dual,
``in-mission'' problem (in contrast with the original ``out-mission'' one),
where the wave is emitted from the source and propagates \textit{backward} in time (being observed in ``usual'' time, it collapses into source):
$t_i = t - \rho (\vect{r}^{(i)};\vect{s})$. In reversed time $T = -t$ these problems are swapped.

If $(\vect{s'};t')$ is a solution of SRP, and $(\vect{s''};t'')$ is a solution of dual problem, then for any $\vect{r}^{(i)}$:

\clin{$t_i = t' + \rho(\vect{r}^{(i)};\vect{s'})$ and $t_i = t'' - \rho(\vect{r}^{(i)};\vect{s''})$, thus
$\rho(\vect{r}^{(i)};\vect{s'}) + \rho(\vect{r}^{(i)};\vect{s''}) = t'' - t' = const$}

\noindent
which may be interpreted as: all sensors belong to the ``ellipsoid'' with $\vect{s'}$ and $\vect{s''}$ being its ``focuses''.
The following example shows that it's possible in $H$.

{\ }

\examp{exampDual}
\begin{mexamp}
Let $H=l_2$, $E = \bigl\{ \vect{x} \in H\mid g(\vect{x}) = \frac{x_1^2}{2} + \sum\limits_{k=2}^{\infty} x_k^2 = 1 \bigr\}$,
and $\vect{s'} = (-1;0;0;...)$, $\vect{s''} = (1;0;0;...)$. We claim that $\forall \vect{x} \in E$: $\| \vect{x} - \vect{s'} \| + \| \vect{x} - \vect{s''} \| = 2\sqrt{2}$.

\begin{mproof}
Consider $L_n = \bigl\{ (x_1; ...; x_n; 0; 0; ...) \mid x_i\in \mathbb{R}, i=\ovline{1,n} \bigr\} \sim \mathbb{R}^n$ and
$E_n = \bigl\{ \vect{x} \in L_n \mid g_n(\vect{x}) = \frac{x_1^2}{2} + \sum\limits_{k=2}^n x_k^2 = 1 \bigr\}$ for $n \greq 2$.
We now show that $\forall \vect{x} \in E_n$: $\| \vect{x} - \vect{s'} \| + \| \vect{x} - \vect{s''} \| = 2\sqrt{2}$, by induction.
For $n = 2$ it holds true since $E_2$, in $L_2$, is the ellipse $x_1^2 / 2 + x_2^2 = 1$ with focuses $(\pm c;0)$, where $c = \sqrt{2 - 1} = 1$.
Suppose it holds true for $n \greq 2$.

For any $\vect{x} \in E_{n+1}$, let $\vect{x}(\vphi) = (x_1; ...; x_{n-1}; x^{(\vphi)}_n; x^{(\vphi)}_{n+1})$,
where $x^{(\vphi)}_n = x_n \cos \vphi - x_{n+1} \sin \vphi$ and $x^{(\vphi)}_{n+1} = x_n \sin \vphi + x_{n+1} \cos \vphi$, for $\vphi \in [0;2\pi)$.
$g_{n+1}(\vect{x}(\vphi)) = \frac{x_1^2}{2} + \sum\limits_{k=2}^{n-1} x_k^2 + (x^{(\vphi)}_n)^2 + (x^{(\vphi)}_{n+1})^2 =
\frac{x_1^2}{2} + \sum\limits_{k=2}^{n+1} x_k^2 = g_{n+1}(\vect{x}) = 1$, hence $\vect{x}(\vphi) \in E_{n+1}$ too.
Similarly, $\| \vect{x}(\vphi) - \vect{s'} \|^2 = \frac{(x_1 + 1)^2}{2} + \sum\limits_{k=2}^{n-1} x_k^2 + (x^{(\vphi)}_n)^2 + (x^{(\vphi)}_{n+1})^2 =
\frac{(x_1 + 1)^2}{2} + \sum\limits_{k=2}^{n+1} x_k^2 = \| \vect{x} - \vect{s'} \|^2$, and $\| \vect{x}(\vphi) - \vect{s''} \| = \| \vect{x} - \vect{s''} \|$.

$\exists \vphi_0 \in [0;2\pi)$: $x^{(\vphi_0)}_{n+1} = 0$, for instance, $\vphi_0 = \begin{cases}
\pi / 2, & x_n = 0,\\
- \arctan \frac{x_{n+1}}{x_n}, & x_n \ne 0,
\end{cases}$ so $\vect{x}(\vphi_0) \in L_n$. Moreover, $g_n(\vect{x}(\vphi_0)) = g_{n+1}(\vect{x}(\vphi_0)) = 1$ $\Rarr$ $\vect{x}(\vphi_0) \in E_n$.
Therefore, $\| \vect{x} - \vect{s'} \| + \| \vect{x} - \vect{s''} \| = \| \vect{x}(\vphi_0) - \vect{s'} \| + \| \vect{x}(\vphi_0) - \vect{s''} \| = 2\sqrt{2}$:
the statement holds true for $n + 1$.

Induction principle leads to it being true for any $n \greq 2$.

Let $\vect{x} \in E$, and consider $\vect{x}^{(n)} = \Bigl( x_1; ...; x_{n-1}; \sqrt{\sum\limits_{k=n}^{\infty} x_k^2}; 0; 0; ... \Bigr) \in L_n$ for $n \greq 3$.
$g_n(\vect{x}^{(n)}) = \frac{x_1^2}{2} + \sum\limits_{k=2}^{n-1} x_k^2 + \sum\limits_{k=n}^{\infty} x_k^2 = g(\vect{x}) = 1$,
thus $\vect{x}^{(n)} \in E_n$ $\Rarr$ $\| \vect{x}^{(n)} - \vect{s'} \| + \| \vect{x}^{(n)} - \vect{s''} \| = 2\sqrt{2}$.

\clin{$\| \vect{x}^{(n)} - \vect{x} \|^2 = \Bigl( \sqrt{\sum\limits_{k=n}^{\infty} x_k^2} - x_n \Bigr)^2 + \sum\limits_{k=n+1}^{\infty} x_k^2
\xrightarrow[n\rarr \infty]{} 0$}

\noindent
and from continuity of metric it follows that

\hfill $2\sqrt{2} = \| \vect{x}^{(n)} - \vect{s'} \| + \| \vect{x}^{(n)} - \vect{s''} \| \xrightarrow[n\rarr \infty]{}
\| \vect{x} - \vect{s'} \| + \| \vect{x} - \vect{s''} \|$ \hfill
\end{mproof}

We place sensors in $E$ as follows: $\vect{r}^{(0)} = (-\sqrt{2};0;0;...)$, $\vect{r}^{(1)} = (\sqrt{2};0;0;...)$, and

\clin{$\vect{r}^{(k)} = (\underbrace{0; ...; 0;}_{k-1} 1; 0; 0; ...)$ for $k \greq 2$}

\noindent
(so $\vect{r}^{(1)} - \vect{r}^{(0)} = (2\sqrt{2};0;0;...)$ and $\vect{r}^{(k)} - \vect{r}^{(0)} = (\sqrt{2}; \underbrace{0; ...; 0;}_{k-2} 1; 0; 0; ...)$ for $k \greq 2$;
$\dot{R} = \{ \hat{\vect{r}}^{(k)} \}_{k\in \mathbb{N}} = \{ \vect{r}^{(k)} - \vect{r}^{(0)} \}_{k\in \mathbb{N}}$ is a basis of $H$).
Since $\forall k\in \mathbb{Z}_+$: $\| \vect{r}^{(k)} - \vect{s'} \| + \| \vect{r}^{(k)} - \vect{s''} \| = 2\sqrt{2}$, we have for $t' = -\sqrt{2}$, $t'' = \sqrt{2}$,
and $t_k = t' + \| \vect{r}^{(k)} - \vect{s'} \|$: $t_k = (t'' - 2\sqrt{2}) + (2\sqrt{2} - \| \vect{r}^{(k)} - \vect{s''} \|) = t'' - \| \vect{r}^{(k)} - \vect{s''} \|$.

In other words, for sensors $R = \{ \vect{r}^{(k)} \}_{k\in \mathbb{Z}_+}$ and moments $\{ t_k \}_{k\in \mathbb{Z}_+}$,
$(\vect{s}';t')$ is the solution of SRP, and $(\vect{s''};t'')$ is the solution of dual problem. \hfill
\end{mexamp}

(It was enough to show that $\forall k\in \mathbb{Z}_+$: $\| \vect{r}^{(k)} - \vect{s'} \| + \| \vect{r}^{(k)} - \vect{s''} \| = 2\sqrt{2}$, without resort to $E$.)

{\ }

Now we return to solving SRP, with the wave propagating \textit{forward} in time.

{\ }

Since $\| \vect{r}^{(i)} - \vect{s} \|^2 = \sum\limits_{j = 1}^{\infty} (r^{(i)}_j - s_j)^2$, $r^{(0)}_j \equiv 0$, and $t_0 = 0$, we arrive to the following:

\clin{$\sum\limits_{j=1}^{\infty} s_j^2 = t^2$, \qquad
$\forall i\in \mathbb{N}$: $\sum\limits_{j=1}^{\infty} \bigl[ (r^{(i)}_j)^2 + s_j^2 - 2 r^{(i)}_j s_j \bigr] = t_i^2 + t^2 - 2 t t_i$}

Subtract the 1st equation from the others, transform and recall that $r^{(i)}_j = 0$ for $j > i$:

\hfill $\begin{cases}
\sum\limits_{j=1}^{\infty} s_j^2 = t^2,\\
\sum\limits_{j=1}^i r^{(i)}_j s_j = \frac{1}{2}\bigl[ \| \vect{r}^{(i)} \|^2 - t_i^2 \bigr] + t t_i, & i\in \mathbb{N}
\end{cases}$ \hfill \equ{impcoordeqset}

Let $b_i = \frac{1}{2}\bigl[ \| \vect{r}^{(i)} \|^2 - t_i^2 \bigr]$, $c_i = t_i$, so $\sum\limits_{j=1}^i r^{(i)}_j s_j = b_i + t c_i$ for all $i\in \mathbb{N}$.

We denote by $A$ the infinite matrix $\| a_{ij} \|_{i, j\in \mathbb{N}} = \| r^{(i)}_j \| = \begin{pmatrix}
r^{(1)}_1 & 0 & 0 & ...\\
r^{(2)}_1 & r^{(2)}_2 & 0 & ...\\
... & ... & ... & \ddots
\end{pmatrix}$.

Also, $S = \begin{pmatrix}
s_1\\
s_2\\
...
\end{pmatrix}$, and $G(t) = B + tC$, where $B = \begin{pmatrix}
b_1\\
b_2\\
...
\end{pmatrix}$, $C = \begin{pmatrix}
c_1\\
c_2\\
...
\end{pmatrix}$. And we have $AS = G(t)$.

The way that we've specified $\{ \vect{r}^{(i)} \}$ allows to express $s_k$ through $t$ from the first $k$ equations of this set;
if we ``cut off'' $A$, $S$, and $G(t)$ after first $k$ rows and columns, the resulting matrix equation $A_k S_k = G_k(t)$
is equivalent to the set of $k$ equations with $k$ unknowns $s_1$, ..., $s_k$.

Using Cramer theorem,

\clin{$s_k = \det \begin{pmatrix}
r^{(1)}_1 & 0 & 0 & ... & 0 & g_1(t)\\
r^{(2)}_1 & r^{(2)}_2 & 0 & ... & 0 & g_2(t)\\
... & ... &... & \ddots &... &...\\
r^{(k)}_1 & r^{(k)}_2 & r^{(k)}_3 & ... & r^{(k)}_{k-1} & g_k(t)\\
\end{pmatrix} / \det A_k =$}

\hfill $= \left\{ \det \begin{pmatrix}
r^{(1)}_1 & 0 & ... & b_1\\
r^{(2)}_1 & r^{(2)}_2 & ... & b_2\\
... & ... & \ddots & ...\\
r^{(k)}_1 & r^{(k)}_2 & ... & b_k\\
\end{pmatrix} + t \det \begin{pmatrix}
r^{(1)}_1 & 0 & ... & c_1\\
r^{(2)}_1 & r^{(2)}_2 & ... & c_2\\
... & ... & \ddots & ...\\
r^{(k)}_1 & r^{(k)}_2 & ... & c_k\\
\end{pmatrix} \right\} / \prod\limits_{i=1}^k r^{(i)}_i = \wtilde{b_k} + t \wtilde{c_k}$ \hfill \equ{coordfromt}

$S = \wtilde{B} + t \wtilde{C}$ with $\wtilde{B} = \begin{pmatrix}
\wtilde{b_1}\\ \wtilde{b_2}\\ ...
\end{pmatrix}$ and $\wtilde{C} = \begin{pmatrix}
\wtilde{c_1} \\ \wtilde{c_2}\\ ...
\end{pmatrix}$; $A(\wtilde{B} + t \wtilde{C}) = B + tC$ $\Rarr$ $A \wtilde{B} = B$, $A \wtilde {C} = C$.

Substituting \equref{coordfromt} into $\sum\limits_{j=1}^{\infty} s_j^2 = t^2$ gives $\sum\limits_{j=1}^{\infty} (\wtilde{b_j} + t \wtilde{c_j})^2 = t^2$ \hfill
\equ{mainquadeq}

{\ }

\textbf{Case 0:} $t = 0$ is a root of \equref{mainquadeq}. We claim that $\vect{s} = \theta$ is the unique solution of SRP then.

\begin{mproof}
$t=0$ turns \equref{mainquadeq} into equality: $\sum\limits_{j=1}^{\infty} \wtilde{b_j}^2 = 0$ $\LRarr$ $\wtilde{b_j} = 0$ for all $j\in \mathbb{N}$.
Since $A \wtilde{B} = B$, it follows that $b_i = 0$ for any $i\in \mathbb{N}$:
$\| \vect{r}^{(i)} \|^2 = t_i^2$ $\LRarr$ $\| \vect{r}^{(i)} \| = |t_i|$.
On the other hand, for any solution $(\vect{s};t)$ of SRP $t_i = t + \| \vect{r}^{(i)} - \vect{s} \| = - \| \vect{s} \| + \| \vect{r}^{(i)} - \vect{s} \|$.
Therefore $\| \vect{r}^{(i)} \| = \bigl| \| \vect{r}^{(i)} - \vect{s} \| - \| \vect{s} \| \bigr|$.

a) $\| \vect{r}^{(i)} \| = \| \vect{r}^{(i)} - \vect{s} \| - \| \vect{s} \|$ $\LRarr$ $\| \vect{r}^{(i)} \| + \| -\vect{s} \| = \| \vect{r}^{(i)} + (- \vect{s}) \|$.

b) $\| \vect{r}^{(i)} \| = - \| \vect{r}^{(i)} - \vect{s} \| + \| \vect{s} \|$ $\LRarr$ $\| \vect{s} \| - \| \vect{r}^{(i)} \| = \| \vect{s} - \vect{r}^{(i)} \|$.

In any case, $\vect{r}^{(i)}$ and $\vect{s}$ are linearly dependent for any $i\in \mathbb{N}$.
Since $\dot{R} = \{ \vect{r}^{(i)} \}_{i\in \mathbb{N}}$ is linearly independent, it is only possible when $\vect{s} = \theta$.
\end{mproof}

{\ }

Until now, the solving method had little relation with infinite dimensionality of $H$.

{\ }

\textbf{Case 1:} $t = 0$ isn't a root of \equref{mainquadeq} (thus $\sum\limits_{j=1}^{\infty} \wtilde{b_j}^2 \ne 0$, and $\vect{s} \ne \theta$).

We divide it by $t$: $\sum\limits_{j = 1}^{\infty} (\wtilde{c_j} + z \wtilde{b_j})^2 = 1$ \hfill \equ{mainquadeqinv}

\noindent
where $z = 1/t < 0$. By assumption, SRP has at least 1 solution, so for some $z_{e;0} = 1/t_{e;0}$ \equref{mainquadeqinv} holds true,
implying $\{ \wtilde{c_j} + z_{e;0} \wtilde{b_j} \}_{j\in \mathbb{N}} = \vect{v} \in H$.

The relations $\wtilde{C} = \vect{v} - z_{e;0} \wtilde{B}$ and $\wtilde{B} = \frac{1}{z_{e;0}} (\vect{v} - \wtilde{C})$
show that $\wtilde{B}$ and $\wtilde{C}$ belong or don't belong to $H$ simultaneously;
the series $\sum\limits_{j=1}^{\infty} \wtilde{b_j}^2$ and $\sum\limits_{j=1}^{\infty} \wtilde{c_j}^2$ both converge or both diverge.

{\ }

\textbf{Subcase 1a} (ruled out in $\mathbb{R}^n$)\textbf{:} $\sum\limits_{j=1}^{\infty} \wtilde{b_j}^2$ and $\sum\limits_{j=1}^{\infty} \wtilde{c_j}^2$ diverge.
Yet for some $z_{e;0}$: $\sum\limits_{j = 1}^{\infty} (\wtilde{c_j} + z_{e;0} \wtilde{b_j})^2$ converges to 1.
Assuming $\exists z'\ne z_{e;0}$ such that $\sum\limits_{j = 1}^{\infty} (\wtilde{c_j} + z' \wtilde{b_j})^2$ converges,
we obtain from equality

\clin{$\wtilde{b_j} = \frac{1}{z_{e;0} - z'} \bigl( (\wtilde{c_j} + z_{e;0} \wtilde{b_j}) - (\wtilde{c_j} + z' \wtilde{b_j}) \bigr)$}

\noindent
the convergence of $\sum\limits_{j=1}^{\infty} \wtilde{b_j}^2$, which contradicts the assumption of the subcase.

Hence $z_{e;0}$ is the one and only value not just satisfying \equref{mainquadeqinv}, but providing the convergence of the series in the left side of \equref{mainquadeqinv}.
How to obtain it? (Recall that we don't allow ourselves to ``go over all $z < 0$ and select the one satisfying the equation'').

$\exists n_0$: $\sum\limits_{j=1}^{n_0} \wtilde{b_j}^2 > 0$, therefore for any $n \greq n_0$: $f_n(z) = \sum\limits_{j=1}^n (\wtilde{c_j} + z \wtilde{b_j})^2 - 1 =$

\clin{$= \bigl[ \sum\limits_{j=1}^n \wtilde{b_j}^2 \bigr] z^2 + \bigl[ 2\sum\limits_{j=1}^n \wtilde{b_j} \wtilde{c_j} \bigr] z + \bigl[ \sum\limits_{j=1}^n \wtilde{c_j}^2 - 1 \bigr] =
\alpha_n z^2 + \beta_n z + \gamma_n$}

\noindent
is a quadratic trinomial with $\alpha_n > 0$. And $f_n(z) \lesq f_{n+1}(z) \lesq f_{\infty}(z)$, so $f_n(z_{e;0}) \lesq 0$: the equation $f_n(z) = 0$ has at least one root.
It is well known that $\{ z \in \mathbb{R} \colon f_n(z) \lesq 0 \}$ is the segment $[z^{(n)}_{-};z^{(n)}_{+}]$ whose center is $z_n = -\frac{\beta_n}{2 \alpha_n}$.

For any $\veps > 0$ the series diverges at $z_{e;0} - \veps$ and $z_{e;0} + \veps$, therefore $\exists n=n(\veps)$: $f_n(z_{e;0} - \veps) > 0$ and $f_n(z_{e;0} + \veps) > 0$.
Consequently, $[z^{(n)}_{-};z^{(n)}_{+}] \subset (z_{e;0} - \veps ; z_{e;0} + \veps)$; in particular, $z_n \in (z_{e;0} - \veps ; z_{e;0} + \veps)$.
That is, $z_n = - \bigl[ \sum\limits_{j=1}^n \wtilde{b_j} \wtilde{c_j} \bigr] / \bigl[ \sum\limits_{j=1}^n \wtilde{b_j}^2 \bigr] \xrightarrow[n\rarr \infty]{} z_{e;0}$.

{\ }

\textbf{Subcase 1b:} $\sum\limits_{j=1}^{\infty} \wtilde{b_j}^2$ and $\sum\limits_{j=1}^{\infty} \wtilde{c_j}^2$ converge.
From the common properties of series it follows that we can rewrite \equref{mainquadeqinv} as

\hfill $\bigl[ \sum\limits_{j=1}^{\infty} \wtilde{b_j}^2 \bigr] z^2 + \bigl[ 2\sum\limits_{j=1}^{\infty} \wtilde{b_j} \wtilde{c_j} \bigr] z + \bigl[ \sum\limits_{j=1}^{\infty} \wtilde{c_j}^2 - 1 \bigr] = 0$ \hfill \equ{mainquadeqinvinf}

\noindent
or $\alpha z^2 + \beta z + \gamma = 0$, where $\alpha > 0$;
what's left to do is to solve it, $z_{\pm} = \frac{-\beta + \sqrt{D}}{2\alpha}$ ($D = \beta^2 - 4 \alpha \gamma \greq 0$ since $z_{e;0}$ is a root),
and select the root(s) $z$ such that $z < 0$ and $t = 1/z \lesq t_i$ for any $i \in \mathbb{N}$.

This concludes the description of the solving method for SRP in $H$.

{\ }

\equref{mainquadeqinvinf} can have 2 distinct roots satisfying $t \lesq t_i \bigl|_{i\in \mathbb{Z_+}}$, --- even when sensors make a basis,
SRP in $H$ can have 2 distinct solutions, as the following example indicates.

{\ }

\examp{exampNonUnique}
\begin{mexamp}
Let non-$\theta$ sensors, $\dot{R} = \{ \vect{r}^{(k)} \}_{k \in \mathbb{N}}$, be $\vect{r}^{(k)} = \frac{1}{k} \vect{e_k}$,
with $\{ \vect{e_k} \}$ being an orthonormal basis of $H$. For the source $\vect{s'} = - \sum\limits_{k=1}^{\infty} \frac{1}{k} \vect{e_k} = (-1;-\frac{1}{2};-\frac{1}{3};...)$,
which emits the wave at the moment $t' = - \| \vect{s'} \| = - \sqrt{\sum\limits_{k=1}^{\infty} \frac{1}{k^2}} = -\frac{\pi}{\sqrt{6}}$,
the moments $\{ t_k \}$ when $k$-th sensor is reached by this wave are such that

\clin{$(t_k - t')^2 = \| \vect{r}^{(k)} - \vect{s'} \|^2 = \sum\limits_{i\in \mathbb{N}, i\ne k} (0 - s'_i)^2 + (\frac{1}{k} - s'_k)^2 =
\sum\limits_{i\in \mathbb{N}} \frac{1}{i^2} + \frac{3}{k^2} = \frac{\pi^2}{6} + \frac{3}{k^2}$}

\noindent
implying $t_k = -\frac{\pi}{\sqrt{6}} + \sqrt{\frac{\pi^2}{6} + \frac{3}{k^2}} = \frac{3}{k^2 \bigl( \frac{\pi}{\sqrt{6}} + \sqrt{\frac{\pi^2}{6} + \frac{3}{k^2}} \bigr)} > 0$
(by construction, for $\vect{r}^{(0)} = \theta$, $t_0 = 0$).

Now we solve the corresponding SRP in accordance with the procedure described above, knowing that $(\vect{s'};t')$ is a solution.
The basis is $\{ \vect{e_i} \}$; the equations from \equref{impcoordeqset} take the form of

\clin{$\frac{1}{i} s_i = \frac{1}{2}\Bigl[ \frac{1}{i^2} - \bigl\{ \frac{3}{i^2 \bigl( \frac{\pi}{\sqrt{6}} + \sqrt{\frac{\pi^2}{6} + \frac{3}{i^2}} \bigr)} \bigr\}^2 \Bigr] +
t\frac{3}{i^2 \bigl( \frac{\pi}{\sqrt{6}} + \sqrt{\frac{\pi^2}{6} + \frac{3}{i^2}} \bigr)}$ $\LRarr$ $s_k = \wtilde{b_k} + t \wtilde{c_k}$}

\noindent
where $\wtilde{b_k} = \frac{1}{k} \cdot\frac{1}{2}\Bigl[ 1 - \frac{9}{k^2 \bigl( \frac{\pi}{\sqrt{6}} + \sqrt{\frac{\pi^2}{6} + \frac{3}{k^2}} \bigr)^2} \Bigr]$,
$\wtilde{c_k}= \frac{1}{k} \cdot\frac{3}{\frac{\pi}{\sqrt{6}} + \sqrt{\frac{\pi^2}{6} + \frac{3}{k^2}}}$.

It is clear that $t = 0$ isn't a root of $\sum\limits_{k=1}^{\infty}(\wtilde{b_k} + t \wtilde{c_k})^2 = t^2$, $\sum\limits_{k=1}^{\infty} \wtilde{b_k}^2$ converges
and $\sum\limits_{k=1}^{\infty} \wtilde{c_k}^2$ converges. Thus we can switch to $z = 1/t$ and the equation $\alpha z^2 + \beta z + \gamma = 0$ from Subcase 1b.
$D \greq 0$ because $z' = 1/t' = -\frac{\sqrt{6}}{\pi}$ is a root. Let $z''$ be a second root; we claim that $z'' < 0$ and $z'' \ne z'$.

\begin{mproof}
$z' z'' = \frac{\gamma}{\alpha}$. Since $\alpha > 0$, $\sign z' z'' = \sign \gamma$. From $\sum\limits_{k=1}^{\infty} \wtilde{c_k}^2 =
\sum\limits_{k=1}^{\infty} \frac{9}{k^2 \bigl( \frac{\pi}{\sqrt{6}} + \sqrt{\frac{\pi^2}{6} + \frac{3}{k^2}} \bigr)^2} >$

\clin{$> \frac{9}{\bigl( \frac{\pi}{\sqrt{6}} + \sqrt{\frac{\pi^2}{6} + 3} \bigr)^2} \sum\limits_{k=1}^{\infty} \frac{1}{k^2} =
\frac{3\pi^2}{2 \bigl( \frac{\pi}{\sqrt{6}} + \sqrt{\frac{\pi^2}{6} + 3} \bigr)^2} = \frac{3}{2 \bigl( \frac{1}{\sqrt{6}} + \sqrt{\frac{1}{6} + \frac{3}{\pi^2}} \bigr)^2} >
\frac{3}{2 \bigl( \frac{1}{\sqrt{6}} + \frac{1}{\sqrt{2}} \bigr)^2} = \frac{9}{(1 + \sqrt{3})^2} > 1$}

\noindent
it follows that $\gamma > 0$, hence $z'' < 0$.

Assume that $z' = z''$, then $\frac{\gamma}{\alpha} = (z')^2 = \frac{6}{\pi^2}$ $\LRarr$
$\sum\limits_{k=1}^{\infty} \wtilde{c_k}^2 = 1 + \frac{6}{\pi^2} \sum\limits_{k=1}^{\infty} \wtilde{b_k}^2$ $\LRarr$
$\sum\limits_{k=1}^{\infty} \bigl( \wtilde{c_k}^2 - \frac{6}{\pi^2} \wtilde{b_k}^2 \bigr) = 1$.

However, when $k \greq 3$, 
$0 < k\wtilde{b_k} \lesq \frac{1}{2}\Bigl[ 1 - \frac{9}{k^2 \bigl( \frac{\pi}{\sqrt{6}} + \sqrt{\frac{\pi^2}{6} + \frac{1}{3}} \bigr)^2} \Bigr] \lesq \frac{1}{2} < 1 <
\frac{3}{\frac{\pi}{\sqrt{6}} + \sqrt{\frac{\pi^2}{6} + \frac{1}{3}}} < k \wtilde{c_k}$ $\Rarr$ $0 < \wtilde{b_k} < \wtilde{c_k}$ $\Rarr$
$0 < \frac{\sqrt{6}}{\pi} \wtilde{b_k} < \wtilde{c_k}$ $\Rarr$ $\wtilde{c_k}^2 > \frac{6}{\pi^2} \wtilde{b_k}^2$.
A little more numerical computation, and we get $\sum\limits_{k=1}^3 \bigl( \wtilde{c_k}^2 - \frac{6}{\pi^2} \wtilde{b_k}^2 \bigr) \approx 1.139918 > 1$,
so $\sum\limits_{k=1}^{\infty} \bigl( \wtilde{c_k}^2 - \frac{6}{\pi^2} \wtilde{b_k}^2 \bigr) > 1$; a contradiction. Therefore $z'' \ne z'$.
\end{mproof}

Then $t'' = 1/z'' < 0 < t_k$, and $\vect{s''} = \{ \wtilde{b_j} + t'' \wtilde{c_j} \}$,
is another solution, different from $(\vect{s'};t')$. \hfill
\end{mexamp}

{\ }

\textbf{Remark.} Non-uniqueness of SRP solution is a well known occasion in $\mathbb{R}^n$.
For example, in $\mathbb{R}^2$ we place 3 sensors on half-hyperbola, and emit the wave at the moment $t'$ from the focus $\vect{s'}$
of hyperbola. Then another focus, $\vect{s''}$, emitting at the moment $t'' = t' + \| \vect{r}^{(i)} - \vect{s'} \| - \| \vect{r}^{(i)} - \vect{s''} \| = t' + const$,
is a different solution of the SRP defined by $\{t_i\}_{i=1}^3$.

{\ }

Below we present some sufficient conditions for uniqueness of SRP solution (Prop. \ref{propDual}--\ref{propExtUnique}).

{\ }

\prop{propDual}
If the dual ``in-mission'' problem, ``$t_i = t - \| \vect{r}^{(i)} - \vect{s} \|$ for any $i\in \mathbb{Z}_+$'',
also has a solution, then the solution of the original SRP is unique.

\begin{mproof}
The implied set of equations \equref{impeqset}, when solved in $t$, has no more than 2 roots, and includes the solution $t' \lesq 0$ of SRP,
along with the solution $t'' \greq 0$ of dual problem. Note that $t'' \ne t'$, otherwise $t' = t'' = 0$, $\vect{s'} = \vect{s''} = \theta$, and the wave reaches $\vect{r}^{(1)}$
at the moment $t_1 \ne 0$ when propagating both forward and backward in time, --- $t_1 = \| \vect{r}^{(1)} \| > 0$ and $t_1 = - \| \vect{r}^{(1)} \| < 0$;
a contradiction. In other words, $t'' > 0$ cannot be a solution of SRP, and $t'$ is the unique solution.
\end{mproof}

{\ }

\prop{propSensor}
If the solution $\vect{s'}$ of the SRP is identical to one of sensors, then it is unique.

\begin{mproof}
Without loss of generality we suppose $\vect{s'} = \vect{r}^{(0)} = \theta$, and $t' = 0$.
Now, assume that $(\vect{s''};t'')$ is another solution. Then for each $\vect{r}^{(i)}$, $i\in \mathbb{N}$, we have $\begin{cases}
t_i = 0 + \| \vect{r}^{(i)} - \theta \|,\\
t_i = t'' + \| \vect{r}^{(i)} - \vect{s''} \|,
\end{cases}$ and for $i = 0$: $0 = t'' + \| \theta - \vect{s''} \|$ $\LRarr$ $t'' = - \| \vect{s''} \|$.
Thus $\| \vect{r}^{(i)} \| = - \| \vect{s''}  \| + \| \vect{r}^{(i)} - \vect{s''} \|$ $\LRarr$ $\| \vect{r}^{(i)} \| + \| - \vect{s''}  \| = \| \vect{r}^{(i)} + (- \vect{s''}) \|$;
linear dep. of $\vect{r}^{(i)}$ and $\vect{s''}$ follows. $\vect{s''} \ne \theta$ then leads to $\vect{r}^{(i)} \in L(\{ \vect{s''} \})$,
for any $i$. This contradicts the linear indep. of $\dot{R}$, so the assumption is wrong.
\end{mproof}

(This conforms with Case 0 above).

{\ }

Of course, we prefer the conditions relating  only to the set of sensors, so that for any position of the source the solution of the SRP is unique and identical to that position.
This is important when sensors must be placed \textit{before} the source appears anywhere in space and emits the wave.

{\ }

\prop{propDivUnique}
If the SRP has a solution, and $\exists \{ n_k \}_{k=1}^{\infty}$, $n_k < n_{k+1}$: $\vect{r}^{(n_k)}\perp \vect{r}^{(i)}$ for $1 \lesq i < n_k$,
and $\| \vect{r}^{(n_k)} \| \in [\lambda; \mu]$ with $\lambda > 0$, then this solution is unique.

\begin{mproof}
We denote the SRP solution by $(\vect{s}';t')$. If $\vect{s'} = \theta = \vect{r}^{(0)}$, it is unique by \propref{propSensor}. Consider $\vect{s'} \ne \theta$.
The basis $B = \{ \vect{e}_k \}$ is made from $\dot{R} = \{ \vect{r}^{(k)} \}$ by Gram-Schmidt orthogonalization:

\clin{$\vect{e}_k = \vect{d}_k / \| \vect{d}_k \|$, where $\vect{d}_k = \vect{r}^{(k)} - \sum\limits_{j=1}^{k-1} {<}\vect{r}^{(k)};\vect{e}_j{>} \vect{e}_j$}

\noindent
hence in $B$, $\vect{r}^{(n_k)}_j = {<}\vect{r}^{(n_k)};\vect{e}_j{>} = 0$ for $j < n_k$, and $r^{(n_k)}_{n_k} = \| \vect{r}^{(n_k)} \|$.

Let $n=n_k$, then by \eqref{coordfromt}: $\wtilde{c_n} = \left|\begin{matrix}
r^{(1)}_1 & 0 & 0 & ... & 0 & t_1\\
r^{(2)}_1 & r^{(2)}_1 & 0 & ... & 0 & t_2\\
... & ... & ... & \ddots & ... & ...\\
r^{(n-1)}_1 & r^{(n-1)}_2 & r^{(n-1)}_3 & ... & r^{(n-1)}_{n-1} & t_{n-1}\\
0 & 0 & 0 & ... & 0 & t_n
\end{matrix}
\right| / \prod\limits_{i=1}^n r^{(i)}_i =$

\noindent
$= t_n (-1)^{n+n} \det A_{n-1} / \prod\limits_{i=1}^n r^{(i)}_i = t_n / r^{(n)}_n = t_n / \| r^{(n)} \|$, so $|\wtilde{c_n}| \greq |t_n| / \mu$.

In turn, $t_n = t' + \| \vect{r}^{(n)} - \vect{s'} \| = -\| \vect{s'} \| + \sqrt{{<}\vect{r}^{(n)} - \vect{s'}; \vect{r}^{(n)} - \vect{s'}{>}} =$

\clin{$= \sqrt{\| \vect{s'} \|^2 + \| \vect{r}^{(n)} \|^2 - 2 {<} \vect{r}^{(n)};\vect{s'}{>}} - \| \vect{s'} \| =
\sqrt{\| \vect{s'} \|^2 + \| \vect{r}^{(n)} \|^2 - 2 \| \vect{r}^{(n)} \| s'_n} - \| \vect{s'} \|$}

$\| \vect{r}^{(n)} \| \lesq \mu$ and $\vect{s'} \in H$ $\Rarr$ $s'_n \xrightarrow[k\rarr \infty]{} 0$,
implying $\| \vect{r}^{(n)} \| s'_n \xrightarrow[k\rarr \infty]{} 0$ ($n = n_k \rarr \infty$ as $k\rarr \infty$).

Therefore $\varliminf\limits_k t_n = \sup\limits_{m\in \mathbb{N}} \inf\limits_{k \greq m} t_n \greq
\bigl[ \sqrt{\| \vect{s'} \|^2 + \lambda^2} - \| \vect{s'} \| \bigr] > 0$, and $\varliminf\limits_k |\wtilde{c_n}| > 0$.
Consequently, $\lim\limits_{k\rarr \infty} \wtilde{c_{n_k}} \ne 0$, so $\sum\limits_{j=1}^{\infty} \wtilde{c_j}^2 = \infty$.
Thus we are in Subcase 1a, where the solution of SRP is unique.
\end{mproof}

The trivial example of such $\dot{R}$ is any \textit{orthonormal} basis of $H$ ($n_k = k$, $\lambda = \mu = 1$).

{\ }

Now, for \textit{arbitrary} basis $\dot{R}$ of $H$, let $\dot{R}'$ be the following ``extension'' of $\dot{R}$:
$\dot{R}' = \dot{R} \cup \{ \vect{r}^{(\omega + 1)} \}$, where $\vect{r}^{(\omega + 1)} = - \vect{r}^{(1)}$. Respectively, $R' = R \cup \{ \vect{r}^{(\omega + 1)} \}$.
The wave reaches this additional sensor, opposite to $\vect{r}^{(1)}$, at the moment $t_{\omega + 1}$.

{\ }

\prop{propExtUnique}
If the SRP defined by $R'$ and $\{ t_i \}_{i\in \mathbb{Z}_+ \cup \{ \omega + 1 \}}$ has a solution, then it is unique.

\begin{mproof}
Assuming the contrary, let $(\vect{s'};t')$ and $(\vect{s''};t'')$ be the distinct solutions of such SRP.
The reasonings above show that $\vect{s}$ is determined uniquely by $t$ ($s_j = \wtilde{b_j} + t \wtilde{c_j}$), therefore $t' \ne t''$.
From \equref{impcoordeqset} for $i = 1$ and $i = \omega + 1$
(it is clear that $\vect{r}^{(\omega + 1)} = (-r^{(1)}_1; 0; 0; ...)$ and $\| \vect{r}^{(\omega + 1)} \| = \| \vect{r}^{(1)} \|$):

\clin{$\begin{cases}
r^{(1)}_1 s_1' = \frac{1}{2} \bigl[ \| \vect{r}^{(1)} \|^2 - t_1^2 \bigr] + t' t_1, & \mathbf{(1')}\\
- r^{(1)}_1 s_1' = \frac{1}{2} \bigl[ \| \vect{r}^{(1)} \|^2 - t_{\omega + 1}^2 \bigr] + t' t_{\omega + 1} & \mathbf{(2')}
\end{cases}$ and $\begin{cases}
r^{(1)}_1 s_1'' = \frac{1}{2} \bigl[ \| \vect{r}^{(1)} \|^2 - t_1^2 \bigr] + t'' t_1, & \mathbf{(1'')}\\
- r^{(1)}_1 s_1'' = \frac{1}{2} \bigl[ \| \vect{r}^{(1)} \|^2 - t_{\omega + 1}^2 \bigr] + t'' t_{\omega + 1} & \mathbf{(2'')}
\end{cases}$}

We subtract $(1'')$ from $(1')$, and $(2'')$ from $(2')$:

\clin{$\begin{cases}
r^{(1)}_1 (s_1' - s_1'') = (t' - t'') t_1,\\
-r^{(1)}_1 (s_1' - s_1'') = (t' - t'') t_{\omega + 1}
\end{cases}$ $\Rarr$ $(t' - t'')(t_1 + t_{\omega + 1}) = 0$ $\Rarr$ $t_{\omega + 1} = - t_1$}

Then we add $(1')$ and $(2')$: $0 = \| \vect{r}^{(1)} \|^2 - t_1^2$ $\LRarr$ $|t_1| = \| \vect{r}^{(1)} \|$. To be definite, suppose $t_1 \greq 0$.
For any solution $(\vect{s};t)$ of the SRP under study (that is, for $(\vect{s'};t')$ and $(\vect{s''};t'')$), we have $t = -\| \vect{s} \|$ and
$\begin{cases}
t_1 = t + \| \vect{r}^{(1)} - \vect{s} \|,\\
t_{\omega + 1} = t + \| \vect{r}^{(\omega + 1)} - \vect{s} \|,
\end{cases}$ hence $\| \vect{r}^{(1)} - (-\vect{r}^{(1)}) \| = 2 t_1 = \| \vect{r}^{(1)} - \vect{s} \| - \| -\vect{r}^{(1)} - \vect{s} \|$. Or:

\clin{$\bigl\| (\vect{r}^{(1)} - \vect{s}) - (-\vect{r}^{(1)} - \vect{s}) \bigr\| = \| \vect{r}^{(1)} - \vect{s} \| - \| -\vect{r}^{(1)} - \vect{s} \|$}

\noindent
which implies linear dependency of $(\vect{r}^{(1)} - \vect{s})$ and $(-\vect{r}^{(1)} - \vect{s})$:
$a (\vect{r}^{(1)} - \vect{s}) + b (-\vect{r}^{(1)} - \vect{s}) = \theta$ $\LRarr$
$(a + b) \vect{s} = (a - b) \vect{r}^{(1)}$.
$a + b \ne 0$, otherwise we divide the 1st equation by $a$ and come to $\vect{r}^{(1)} - \vect{s} + \vect{r}^{(1)} + \vect{s} = \theta$,
or $\vect{r}^{(1)} = \theta$, --- a contradiction. Thus we can divide the 2nd equation by $(a + b)$:
$\vect{s} = \frac{a - b}{a + b} \vect{r}^{(1)} \in L(\{ \vect{r}^{(1)} \})$, therefore $\vect{s} = (s_1; 0; 0; ...)$,
and $\| \vect{r}^{(1)} - \vect{s} \| = |r^{(1)}_1 - s_1|$, $\| \vect{r}^{(\omega + 1)} - \vect{s} \| = |r^{(1)}_1 + s_1|$.

$t_1 \greq 0$ $\Rarr$ $|r^{(1)}_1 + s_1| \lesq |r^{(1)}_1 - s_1|$.
The basis $B = \{ \vect{e}_i \}$ was made from $\dot{R}$ by Gram-Schmidt orthogonalization, hence $r^{(1)}_1 > 0$, and $r^{(\omega + 1)}_1 = -r^{(1)}_1 < 0$.
Thus $s_1 \lesq 0$ $\Rarr$ $t = -\| \vect{s} \| = s_1$.

We now take into account other sensors; one is enough, for instance, $\vect{r}^{(2)} = (r^{(2)}_1; r^{(2)}_2; 0; ...) = (p; h; 0; ...)$, where $h \ne 0$.
$t_2 = t + \| \vect{r}^{(2)} - \vect{s} \| = t + \sqrt{(p - s_1)^2 + h^2}$ turns into equality for $(\vect{s'};t')$ and $(\vect{s''};t'')$:
$t' + \sqrt{(p - t')^2 + h^2} = t'' + \sqrt{(p - t'')^2 + h^2}$.

However, $f(t) = t + \sqrt{(t - p)^2 + h^2}$ has the derivative $f'(t) = 1 + \frac{t - p}{\sqrt{(t - p)^2 + h^2}} > 0$, because
$|t - p| < \sqrt{(t - p)^2 + h^2}$. Hence $f(t)$ is strictly increasing, and it must be $t' = t''$, --- a contradiction.
Consequently, the initial assumption is wrong; the solution is unique.
\end{mproof}

{\ }

\prop{propSolNotBoth}
If the SRP \equref{emeqset} has a solution, and $(\vect{s''};t'')$ is a different solution of implied \equref{impeqset}, 
then $(\vect{s''};t'')$ is either the solution of SRP, or the solution of dual problem.

\begin{mproof}
Assume the contrary, then $\exists m,k \in \mathbb{Z}_+$: $\begin{cases}
t_m = t'' - \| \vect{r}^{(m)} - \vect{s''} \|,\\
t_k = t'' + \| \vect{r}^{(k)} - \vect{s''} \|.
\end{cases}$
Let $(\vect{s'};t')$ be the solution of SRP \equref{emeqset}, so $\forall i\in \mathbb{Z}_+$: $t_i = t' + \| \vect{r}^{(i)} - \vect{s'} \|$.
In particular, $\begin{cases}
t_m = t' + \| \vect{r}^{(m)} - \vect{s'} \|,\\
t_k = t' + \| \vect{r}^{(k)} - \vect{s'} \|.
\end{cases}$ Therefore
$\| \vect{r}^{(m)} - \vect{s'} \| + \| \vect{r}^{(m)} - \vect{s''} \| = t'' - t' = \| \vect{r}^{(k)} - \vect{s'} \| - \| \vect{r}^{(k)} - \vect{s''} \|$.

By triangle inequality, $\| \vect{r}^{(m)} - \vect{s'} \| + \| \vect{r}^{(m)} - \vect{s''} \| \greq \| \vect{s'} - \vect{s''} \|$;
on the other hand,

\clin{$\bigl| \| \vect{r}^{(k)} - \vect{s'} \| - \| \vect{r}^{(k)} - \vect{s''} \| \bigr| \lesq \| \vect{s'} - \vect{s''} \|$}

Hence $\begin{cases}
\| \vect{s'} - \vect{r}^{(m)} \| + \| \vect{r}^{(m)} - \vect{s''} \| = \| \vect{s'} - \vect{s''} \|,\\
\| \vect{r}^{(k)} - \vect{s'} \| - \| \vect{r}^{(k)} - \vect{s''} \| = \| \vect{s''} - \vect{s'} \|.
\end{cases}$ From the 1st equality we obtain linear dependency of $\vect{s'} - \vect{r}^{(m)}$ and $\vect{r}^{(m)} - \vect{s''}$:
$\exists a, b$: $a (\vect{s'} - \vect{r}^{(m)}) + b (\vect{r}^{(m)} - \vect{s''}) = \theta$ $\LRarr$ $(b-a)\vect{r}^{(m)} = b \vect{s''} - a \vect{s'}$.
$a\ne b$, otherwise we could divide by $a$ and get $\vect{s'} - \vect{s''} = \theta$; so $\vect{r}^{(m)} = \frac{1}{b-a} (b \vect{s''} - a \vect{s'}) \in L(\{ \vect{s'};\vect{s''} \})$.

From 2nd equality: $\vect{r}^{(k)} - \vect{s'}$ and $\vect{r}^{(k)} - \vect{s''}$ are linearly dependent, $a (\vect{r}^{(k)} - \vect{s'}) + b (\vect{r}^{(k)} - \vect{s''}) = \theta$
$\LRarr$ $(a+b)\vect{r}^{(k)} = a \vect{s'} + b \vect{s''}$, $a + b \ne 0$ or it would be $\vect{s'} - \vect{s''} = \theta$, thus $\vect{r}^{(k)}\in L(\{ \vect{s'};\vect{s''} \})$.

Moreover, for any $j\in \mathbb{Z}_+$ such that $t_j = t'' + \| \vect{r}^{(j)} - \vect{s''} \|$ we can repeat these reasonings for the same $m$, but taking $j$ instead of $k$.
Consequently, $R_+ = \bigl\{ \vect{r}^{(j)} \mid t_j = t'' + \| \vect{r}^{(j)} - \vect{s''} \| \bigr\} \subseteq L(\{ \vect{s'}; \vect{s''} \})$.
Similarly, keeping $k$ and going over suitable $m$, $R_- = \{ \vect{r}^{(j)} \mid t_j = t'' - \| \vect{r}^{(j)} - \vect{s''} \| \} \subseteq L(\{ \vect{s'}; \vect{s''} \})$.
Since $R = R_+ \cup R_-$, we have $R \subseteq L(\{ \vect{s'}; \vect{s''} \})$, which is impossible, because $\dot{R} = R \backslash \{ \theta \}$ is a basis of $H$,
while $\dim L(\{ \vect{s'}; \vect{s''} \}) \lesq 2$. This contradiction proves that the assumption is wrong.
\end{mproof}

In other words, when a solution of SRP exists, the transition from \equref{emeqset} to \equref{impeqset} may add only the solution of dual problem, not some ``mixed'' one.

{\ }

When we have the countable set $R$ of sensors and corresponding moments $\{ t_i \}_{i\in \mathbb{Z}}$, we may ``downdimension'' the original SRP
by taking into account only the sensors from 0-th to $n$-th, $R_n = \{ \vect{r}^{(i)} \}_{i=0}^n$.
Since $\vect{r}^{(0)} = \theta$, we have $L_n := L(R_n) = L(\dot{R}_n)$ and is isomorphic to $\mathbb{R}^n$.

Further, we seek the solution $(\vect{s};t)$ of the problem ``$t_i = t + \| \vect{r}^{(i)} - \vect{s} \|$ for any $i=\ovline{0,n}$'' inside $L_n$.
We denote this ``downdimensioned'' problem by SRP$_n$.

{\ }

\prop{propDownDim}
If the SRP in $H$ has a solution, then SRP$_n$ has a solution for any $n\in \mathbb{N}$.

\begin{mproof}
Denote the solution of original SRP by $(\vect{s}^{(\infty)};t^{(\infty)})$. By projection theorem, $\vect{s}^{(\infty)} = \vect{u} + \vect{h}$,
where $\vect{u} \in L_n$ and $\vect{h} \perp L_n$. If $\vect{h} = \theta$, then $\vect{s}^{(\infty)}$ is the solution of SRP$_n$.
We consider another case, $\vect{h} \ne \theta$. Let $h = \| \vect{h} \|$.

Let $L_n' = L_n \oplus L(\{ \vect{h} \}) = L(R_n \cup \{ \vect{h} \})$. It is isomorphic to $\mathbb{R}^{n+1}$, thus $\vect{x} \in L_n'$
may be written as $(x_1;...;x_n;x_{n+1})$ in the basis made, using Gram-Schmidt orthogonalization, from $\dot{R}_n \cup \{ \vect{h} \}$.
In particular, for $i=\ovline{1,n}$ the sensor $\vect{r}^{(i)}$ has the coordinates $\{ r^{(i)}_j \}_j$ with $r^{(i)}_{n+1} = 0$.
Also, $\vect{s}^{(\infty)} \in L_n'$ and $\vect{s}^{(\infty)} = (s^{(\infty)}_1;...;s^{(\infty)}_n;s^{(\infty)}_{n+1})$ with $s^{(\infty)}_{n+1} = h$.

We now consider the SRP defined by $(R_n;\{ t_i \}_{i=0}^n)$ in $L_n'$; it has (at least one) solution $(\vect{s}^{(\infty)};t^{(\infty)})$.
Following the way of \equref{emeqset}-\equref{impeqset}-\equref{impcoordeqset}-\equref{coordfromt}-\equref{mainquadeq}
(now there's a finite sum instead of series),

\clin{$s_j = \wtilde{b_j} + t \wtilde{c_j}$ for $j=\ovline{1,n}$; $\sum\limits_{j=1}^n (\wtilde{b_j} + t \wtilde{c_j})^2 + s_{n+1}^2 = t^2$}

We rewrite the latter equation, in $t$, as $\alpha t^2 + \beta t + \gamma + s_{n+1}^2 = 0$. Note that $\gamma = \sum\limits_{j=1}^n \wtilde{b_j}^2 \greq 0$.
This equation has the solution $t = t^{(\infty)} \lesq t_i$, $i=\ovline{0,n}$, when $s_{n+1} = \pm h$
(so, actually, this SRP has at least 2 solutions in $L_n'$, symmetrical with respect to $L_n$).

We claim that it has a solution $t^{(n)} \lesq t_i$ when $s_{n+1} = 0$. Consider the cases:

\textbf{Case $\alpha > 0$:} $f_{h^2} (t) = \alpha t^2 + \beta t + \gamma + h^2 = 0$ has a root $t^{(\infty)} \lesq t_i$; if $t'$ is the lesser root of this equation, then all the more
$t' \lesq t_i$. Since $f_{h^2}(t)$ is a quadratic trinomial, for $h^2$ replaced by $0$ it has 2 roots, with the lesser one $t'' < t'$. Let $t^{(n)} = t''$.

\textbf{Case $\alpha < 0$:} $f_0(t)$ has a root(s) because $D = \beta^2 - 4 \alpha (\gamma + 0) \greq \beta^2 \greq 0$ (perhaps the root is multiple).
Its roots are $t_{\pm} = \frac{-\beta \pm \sqrt{D}}{2\alpha}$, and $t_+ t_- = \frac{\gamma}{\alpha} \lesq 0$, thus $t_+ \lesq 0$ (and $t_- \greq 0$).
In other words, there's only 1 root satisfying $t \lesq 0$, which distinguishes the solution of SRP from the solution of dual, ``in-mission'' problem.

Now, if we repeat the solving method after re-enumerating the sensors so that $i$-th sensor ($i=\ovline{1,n}$) becomes $\vect{r}^{(0)}$,
and moving ``the origin of space and time'' to this new $\vect{r}^{(0)}$, then we come to essentially the same SRP, $L_n$, $L_n'$, ...
in different reference frame. And we obtain the single root $T_+ \lesq 0$. But \textit{$T_+$ and $t_+$ are the same moment} of time,
only in different temporal reference frames. Therefore, $T_+ \lesq 0$ means $t_+ \lesq t_i$. Let $t^{(n)} = t_+$.

\textbf{Case $\alpha = 0$, $\beta \ne 0$:} $f_{h^2}(t) = \beta t + \gamma + h^2 = 0$ $\LRarr$ $t = t^{(\infty)} = -\frac{\gamma + h^2}{\beta}$.
$t^{(\infty)} \lesq 0$ $\Rarr$ $\beta > 0$, thus $f_0(\hat{t}) = 0$ for $\hat{t} = -\frac{\gamma}{\beta} \lesq 0$.
Similarly, the symmetry implies $\hat{t} \lesq t_i$ for any $i=\ovline{1,n}$. Let $t^{(n)} = \hat{t}$.

\textbf{Case $\alpha = 0$, $\beta = 0$:} impossible, because $\gamma + h^2 > 0$.

Anyway, $\exists t^{(n)} \lesq t_i$ for any $i=\ovline{0,n}$: $\sum\limits_{j=1}^n (\wtilde{b_j} + t^{(n)} \wtilde{c_j})^2 = (t^{(n)})^2$.
It determines the solution $\vect{s}^{(n)} = \{ \wtilde{b_j} + t^{(n)} \wtilde{c_j} \}_{j=1}^n \in L_n$ of SRP$_n$.
\end{mproof}

{\ }

The statement of \propref{propDownDim} remains true if we take arbitrary finite $\widehat{R} = \{ \vect{r}^{(i_1)};...;\vect{r}^{(i_n)} \} \subset \dot{R}$
and seek the solution of the ``truncated'' SRP ``$t_{i_j} = t + \| \vect{r}^{(i_j)} - \vect{s} \|$ for $j=\ovline{0,n}$'' in $L(\widehat{R})$, ---
just re-enumerate elements of $\dot{R}$ so that $\widehat{R} = \{ \vect{r}^{(1)};...;\vect{r}^{(n)} \}$,
$\dot{R}\backslash \widehat{R} = \{ \vect{r}^{(n+1)}; \vect{r}^{(n+2)}; ... \}$ to get SRP$_n$.

However, this statement becomes false for infinite $\widehat{R} \subset \dot{R}$, in general case. Consider

{\ }

\examp{exampDownDimInfinite}
\begin{mexamp}
Let $\dot{R} = B$ be an orthonormal basis of $H$, $\vect{r}^{(i)} = \vect{e}_i$, thus $r^{(i)}_j = \delta_{ij}$;
also, let $\vect{s'} = \vect{r}^{(1)} = (1;0;0;...)$ and $t' = -1$. Then $t_0 = t' + \| \vect{s'} \| = 0$, $t_1 = t' + \| \vect{r}^{(1)} - \vect{s'} \| = -1$,
and $\forall i \greq 2$: $t_i = t' + \| \vect{r}^{(i)} - \vect{s'} \| = -1 + \sqrt{2}$.

Obviously, $(\vect{s'};t')$ is the solution of the SRP defined by $\bigl(R;\{t_i \}_{i\in \mathbb{Z}_+}\bigr)$.
We claim that for any infinite $\widehat{R} \subset \dot{R}$
such that $\vect{r}^{(1)} \notin \widehat{R}$ the truncated SRP ``$t_i = t + \| \vect{r}^{(i)} - \vect{s} \|$ for any $\vect{r}^{(i)} \in \widehat{R} \cup \{ \vect{r}^{(0)} \}$''
has no solution in $\widehat{L} = L(\widehat{R})$.

\begin{mproof}
Assume the contrary and enumerate the elements of $\widehat{R}$ as the subsequence of $\dot{R}$, ascending:
$\widehat{R} = \{ \hat{\vect{r}}^{(1)}; \hat{\vect{r}}^{(2)}; ... \}$.
$\widehat{R}$ is the orthonormal basis of $\widehat{L}$, in itself $\hat{r}^{(i)}_j = \delta_{ij}$ as well, and $\dim \widehat{L} = \infty$.
We denote the solution of truncated SRP in $\widehat{L}$ by $(\vect{s};t)$, with $\vect{s} = (s_1; s_2; ...)$.

\equref{emeqset}-\equref{impeqset}-\equref{impcoordeqset}-\equref{coordfromt}-\equref{mainquadeq} implies

\clin{$\sum\limits_{k=1}^{\infty} s_k^2 = t^2$, \qquad $\forall k\in \mathbb{N}$:
$\sum\limits_{j=1}^k \hat{r}^{(k)}_j s_j = s_k = \frac{1}{2}(\| \hat{\vect{r}}^{(k)} \|^2 - \hat{t}_k^2) + t \hat{t}_k$}

\noindent
(that is, $b_k = \wtilde{b_k}$, $c_k = \wtilde{c_k}$).
Hence $s_k \equiv \frac{1}{2}\bigl(1 - (\sqrt{2} - 1)^2\bigr) + t(\sqrt{2} - 1) = (\sqrt{2} - 1)(t + 1)$,
therefore $\vect{s}\in H$ only if $s_k \equiv 0$ $\LRarr$ $t = -1$. Then for $\hat{\vect{r}}^{(0)} = \theta = \vect{s}$:
$t_0 = t = -1 \ne 0$, --- a contradiction.
\end{mproof}
\end{mexamp}

{\ }

\prop{propDownDimConv}
If the solution $(\vect{s}^{(\infty)};t^{(\infty)})$, $\vect{s}^{(\infty)} \ne \theta$, of SRP is unique, and for each $n\in \mathbb{N}$
the solution $(\vect{s}^{(n)};t^{(n)})$, $\vect{s}^{(n)} \ne \theta$, of SRP$_n$ is unique, and $\sum\limits_{j=1}^{\infty} \wtilde{c_j}^2 < \infty$, then

\clin{$t^{(n)} \xrightarrow[n\rarr \infty]{} t^{(\infty)}$ and $\vect{s}^{(n)} \xrightarrow[n\rarr \infty]{} \vect{s}^{(\infty)}$}

\begin{mproof}
\equref{coordfromt} gives $\vect{s}^{(\infty)} = (s^{(\infty)}_1; s^{(\infty)}_2; ...) =
( \wtilde{b_1} + t^{(\infty)} \wtilde{c_1}; ...; \wtilde{b_n} + t^{(\infty)} \wtilde{c_n}; s^{(\infty)}_{n+1}; s^{(\infty)}_{n+2}; ...)$
and $\vect{s}^{(n)} = (\wtilde{b_1} + t^{(n)} \wtilde{c_1}; ...; \wtilde{b_n} + t^{(n)} \wtilde{c_n}; 0; 0; ...)$, consequently

\clin{$\| \vect{s}^{(n)} - \vect{s}^{(\infty)} \|^2 = (t^{(n)} - t^{(\infty)})^2 \sum\limits_{j=1}^n \wtilde{c_j}^2 + \sum\limits_{j=n+1}^{\infty} (s^{(\infty)}_j)^2$}

$\sum\limits_{j=n+1}^{\infty} (s^{(\infty)}_j)^2 \xrightarrow[n\rarr \infty]{} 0$ and
$\sum\limits_{j=1}^n \wtilde{c_j}^2 \xrightarrow[n\rarr \infty]{} \sum\limits_{j=1}^{\infty} \wtilde{c_j}^2$; it remains to prove that

\clin{$t^{(n)} \xrightarrow[n\rarr \infty]{} t^{(\infty)}$ $\LRarr$ $z^{(n)} \xrightarrow[n\rarr \infty]{} z^{(\infty)}$}

\noindent
where $z^{(\infty)} = 1/t^{(\infty)}$, $z^{(n)} = 1/t^{(n)}$.

From the assumptions of this proposition it follows that, speaking of SRP, we're in Subcase 1b,
where $z^{(\infty)}$ is one of two roots, $z^{(\infty)}_{\pm} = \frac{-\beta \pm \sqrt{\beta^2 - 4 \alpha \gamma}}{2 \alpha}$, of \equref{mainquadeqinvinf}.
Now, using the symbols $\alpha_n$, $\beta_n$ and $\gamma_n$ from Subcase 1a (this is different from notation used while proving \propref{propDownDim}),
we state that, similarly, $z^{(n)}$ is one of two roots, $z^{(n)}_{\pm} = \frac{-\beta_n \pm \sqrt{\beta_n^2 - 4 \alpha_n \gamma_n}}{2 \alpha_n}$,
of the equation $\alpha_n z^2 + \beta_n z + \gamma_n = 0$, which appears while solving SRP$_n$.

$\alpha_n \rarr \alpha > 0$, $\beta_n \rarr \beta$, $\gamma_n \rarr \gamma$ as $n \rarr \infty$,
therefore $z^{(n)}_- \rarr z^{(\infty)}_-$, $z^{(n)}_+ \rarr z^{(\infty)}_+$, and the selection of the root $z^{(n)}$ in SRP$_n$ ($z^{(n)}_-$ or $z^{(n)}_+$)
becomes the same as the selection of the root $z^{(\infty)}$ in SRP (perhaps for $n \greq n_0$). In any case, $z^{(n)} \rarr z^{(\infty)}$ as $n\rarr \infty$.
\end{mproof}

This proposition shows another method, one of Galerkin kind, to obtain the SRP solution.

\section{SR on unit sphere in Hilbert space}

\forceparindent
Let it be $S = \{ \vect{x} \in H\colon \| \vect{x}\| = 1 \}$. Instead of ``embracing-space-induced'' $\| \vect{x} - \vect{y} \|$,
we consider the so-called geodesic metric, $d \colon S \times S \rarr \mathbb{R}_+$: $d(\vect{x}; \vect{y}) = \arccos {<}\vect{x}; \vect{y}{>} \in [0;\pi]$.

Which (we remind) is really a metric.
\begin{mproof}
Obviously, $d(\vect{x};\vect{x}) = \arccos 1 = 0$ and $d (\vect{x};\vect{y}) = d (\vect{y};\vect{x})$.
If $d (\vect{x};\vect{y}) = 0$, then ${<}\vect{x};\vect{y}{>} = 1 = \| \vect{x} \| \cdot \| \vect{y} \|$ $\Rarr$ $\vect{x}$ and $\vect{y}$ are linearly dependent
with $\vect{y} = a \vect{x}$, $a \greq 0$; $1 = {<}\vect{x};\vect{y}{>} = a \| \vect{x} \|^2 = a$ $\Rarr$  $\vect{x} = \vect{y}$.

The triangle inequality $\forall \vect{x}, \vect{y}, \vect{z} \in S$: $d(\vect{x};\vect{z}) \lesq d(\vect{x};\vect{y}) + d(\vect{y};\vect{z})$
can be established as follows. It is equivalent to $\left[ \begin{array}{l}
d(\vect{x};\vect{y}) + d(\vect{y};\vect{z}) \greq \pi,\\
\cos d (\vect{x};\vect{z}) \greq \cos \bigl( d(\vect{x};\vect{y}) + d(\vect{y};\vect{z}) \bigr), \quad d(\vect{x};\vect{y}) + d(\vect{y};\vect{z}) \lesq \pi;
\end{array}\right.$ we rewrite the inequality in the 2nd case as

\clin{${<}\vect{x};\vect{z}{>} \greq {<}\vect{x};\vect{y}{>} {<}\vect{y};\vect{z}{>} - \sqrt{1 - {<}\vect{x};\vect{y}{>}^2} \cdot \sqrt{1 - {<}\vect{y};\vect{z}{>}^2}$ $\LRarr$}

\clin{$\LRarr$ $\left[ \begin{array}{l}
{<}\vect{x};\vect{y}{>} {<}\vect{y};\vect{z}{>} \lesq {<}\vect{x};\vect{z}{>},\\
\bigl( 1 - {<}\vect{x};\vect{y}{>}^2 \bigr) \bigl( 1 - {<}\vect{y};\vect{z}{>}^2 \bigr) \greq \bigl( {<}\vect{x};\vect{y}{>} {<}\vect{y};\vect{z}{>} - {<}\vect{x};\vect{z}{>} \bigr)^2, \quad {<}\vect{x};\vect{y}{>} {<}\vect{y};\vect{z}{>} \greq {<}\vect{x};\vect{z}{>}
\end{array} \right.$}

The 2nd inequality here, being rearranged,

\hfill $1 + 2 {<}\vect{x};\vect{y}{>} {<}\vect{y};\vect{z}{>} {<}\vect{x};\vect{z}{>} \greq {<}\vect{x};\vect{y}{>}^2 + {<}\vect{y};\vect{z}{>}^2 + {<}\vect{x};\vect{z}{>}^2$ \hfill \equ{triscalarineq}

Using projection theorem (and $\dim H = \infty$), we represent $\vect{y} = y_1 \vect{x} + y_2 \vect{h}$, where $\vect{h} \in S$, $\vect{h} \perp \vect{x}$,
and $\vect{z} = z_1 \vect{x} + z_2 \vect{h} + z_3 \vect{w}$, where $\vect{w} \in S$, $\vect{w} \perp \vect{x}$, $\vect{w} \perp \vect{h}$.
And $y_1^2 + y_2^2 = 1$, $z_1^2 + z_2^2 + z_3^2 = 1$, so
\equref{triscalarineq} $\LRarr$ $1 + 2 y_1 (y_1 z_1 + y_2 z_2) z_1 \greq y_1^2 + (y_1 z_1 + y_2 z_2)^2 + z_1^2$ $\LRarr$
$1 - y_1^2 \greq  z_1^2 (1 - y_1^2) + y_2^2 z_2^2$ $\LRarr$ $y_2^2 z_3^2 \greq 0$
\end{mproof}

{\ }

The set of sensors $R = \{ \vect{r}^{(i)} \}_i \subset S$ (obviously, $\theta \notin R$).
As before, we assume the existence of at least one solution $(\vect{s}_0;t_{e;0})$, $\vect{s}_0 \in S$,
of the SRP ``$t_i = t + d(\vect{r}^{(i)};\vect{s})$ for any $i$''.

{\ }

\textbf{Remark.} We may consider the wave to ``oscillate forever'' on $S$, from $\vect{s}_0$ to antipodal $-\vect{s}_0$ ($d(\vect{s}_0;-\vect{s}_0) = \pi$),
then back to $\vect{s}_0$, and so forth. Then $t_i$ is the \textit{first} time when the wave reaches $\vect{r}^{(i)}$.
However, the wave as the sphere of increasing radius $t - t_{e;0}$ vanishes at $-\vect{s}_0$.

{\ }

The reasonings we've used for the entire $H$ show that if $L(R) \ne H$ and $\vect{s}_0\notin L(R)$,
then the solution is certainly not unique: for $\vect{s}_0 = \vect{u_0} + \vect{u_1}$ with $\vect{u_0} \in L(R)$, $\vect{u_1} \perp L(R)$,
the ``$\vect{s}(\vphi) = \vect{u_0} + \| \vect{u_1} \| (\cos \vphi \cdot \wtilde{\vect{u_1}} + \sin \vphi \cdot \wtilde{\vect{u_2}})$''
(where $\wtilde{\vect{u_2}} \perp L(R), \vect{u_1}$) construction works as well, since $\vect{s}(\vphi) \in S$
and $d(\vect{r}^{(i)};\vect{s}(\vphi)) = \arccos\bigl[ {<}\vect{r}^{(i)};\vect{u_0}{>} + \| \vect{u_1} \| \cos \vphi {<}\vect{r}^{(i)};\wtilde{\vect{u_1}}{>} + 
\| \vect{u_1} \| \sin \vphi {<}\vect{r}^{(i)};\wtilde{\vect{u_2}}{>}\bigr] = \arccos {<}\vect{r}^{(i)};\vect{u_0}{>} =
\arccos {<}\vect{r}^{(i)};\vect{u_0} + \vect{u_1}{>} = d(\vect{r}^{(i)};\vect{s}_0)$.

Therefore, let $R = \{ \vect{r}^{(i)} \}_{i\in \mathbb{N}}$ be a basis of $H$,
and let $B$ be the orthonormal basis of $H$, derived from $R$ by Gram-Schmidt orthogonalization;
thus, in $B$, $\vect{r}^{(i)} = (r^{(i)}_1; ...; r^{(i)}_i; 0; 0; ...)$ (and $r^{(1)}_1 = 1$).

Then $\vect{s}$ can be written in the form of $(s_1;s_2;...)$.

{\ }

Since $t \lesq t_i$ and $t_i - t = d(\vect{r}^{(i)};\vect{s}) \lesq \pi$, we have $t \in [\sup \{ t_i \} - \pi; \inf \{ t_i \}] = \mathbf{\Delta}$ ($|\mathbf{\Delta}| \lesq \pi$).
The equations of SRP are equivalent to $\cos (t_i - t) = {<}\vect{r}^{(i)};\vect{s}{>}$. Adding ``$\vect{s} \in S$'', we have

\hfill
$\sum\limits_{j=1}^{\infty} s_j^2 = 1$,\qquad
$\forall i\in \mathbb{N}$: $\sum\limits_{j=1}^i r^{(i)}_j s_j = \cos t \cos t_i + \sin t \sin t_i$
\hfill \equ{sphcoordeqset}

Similarly to how we've moved from \equref{impcoordeqset} to \equref{coordfromt} and \equref{mainquadeq},
we obtain $s_j = \wtilde{p_j} \cos t + \wtilde{q_j} \sin t$, where

\hfill $\wtilde{p_k} = \left| \begin{matrix}
r^{(1)}_1 & 0 & ... & \cos t_1\\
r^{(2)}_1 & r^{(2)}_2 & ... & \cos t_2\\
... & ... & \ddots & ...\\
r^{(k)}_1 & r^{(k)}_2 & ... & \cos t_k\\
\end{matrix} \right| / \prod\limits_{i=1}^k r^{(i)}_i$, $\wtilde{q_k} = \left| \begin{matrix}
r^{(1)}_1 & 0 & ... & \sin t_1\\
r^{(2)}_1 & r^{(2)}_2 & ... & \sin t_2\\
... & ... & \ddots & ...\\
r^{(k)}_1 & r^{(k)}_2 & ... & \sin t_k\\
\end{matrix} \right| / \prod\limits_{i=1}^k r^{(i)}_i$ \hfill \equ{sphcoeffcoordfromt}

\noindent
and $\sum\limits_{j=1}^{\infty} (\wtilde{p_j} \cos t + \wtilde{q_j} \sin t)^2 = 1$ \hfill \equ{sphmainquadeq}

\textbf{Case 1a:} $\sum\limits_{j=1}^{\infty} \wtilde{p_j}^2$ converges, $\sum\limits_{j=1}^{\infty} \wtilde{q_j}^2$ diverges.
If the series in the left side of \equref{sphmainquadeq} converges for $t$
such that $\sin t \ne 0$, then $\{ \wtilde{q_j} \} = \frac{1}{\sin t} \bigl( \{ \wtilde{p_j} \cos t + \wtilde{q_j} \sin t \} - \{ \wtilde{p_j} \cos t \} \bigr) \in H$,
which contradicts the assumption. Thus, if $t$ satisfies \equref{sphmainquadeq}, then $\sin t = 0$
(in particular, $\sin t_{e;0} = 0$, so $\sum\limits_{j=1}^{\infty} \wtilde{p_j}^2 = 1$).

We go over $t = \pi m \in \mathbf{\Delta}$, $m\in \mathbb{Z}$ (there's 1 or 2 such values),
and select those satisfying \equref{sphmainquadeq} (in fact, they all do).
Then $s_j = \wtilde{p_j} \cos t + \wtilde{q_j} \sin t = \wtilde{p_j} \cos t$ for all $j\in \mathbb{N}$.

\textbf{Case 1b:} $\sum\limits_{j=1}^{\infty} \wtilde{p_j}^2$ diverges, $\sum\limits_{j=1}^{\infty} \wtilde{q_j}^2$ converges.
Similarly, the convergence of the series in the left side of \equref{sphmainquadeq} leads to $\cos t = 0$ ($\sum\limits_{j=1}^{\infty} \wtilde{q_j}^2 = 1$ since it converges
when $t = t_{e;0}$), and we take $t \in  (\frac{\pi}{2} + \pi \mathbb{Z}) \cap \mathbf{\Delta}$, satisfying \equref{sphmainquadeq}.

\textbf{Case 2:} $\sum\limits_{j=1}^{\infty} \wtilde{p_j}^2$ and $\sum\limits_{j=1}^{\infty} \wtilde{q_j}^2$ diverge.
We know that \equref{sphmainquadeq} holds true for $t' = t_{e;0}$; if the series in the left side converges for $t'' \in \mathbf{\Delta}$, $t'' \ne t'$, then we have

\clin{$\begin{cases}
\cos t' \{ \wtilde{p_j} \} + \sin t' \{ \wtilde{q_j} \} = \vect{v'} \in S,\\
\cos t'' \{ \wtilde{p_j} \} + \sin t'' \{ \wtilde{q_j} \} = \vect{v''} \in H,
\end{cases}$ $\Rarr$ $\begin{cases}
\cos t' \sin t'' \{ \wtilde{p_j} \} + \sin t' \sin t'' \{ \wtilde{q_j} \} = \sin t'' \vect{v'} \in H,\\
\cos t'' \sin t' \{ \wtilde{p_j} \} + \sin t'' \sin t' \{ \wtilde{q_j} \} = \sin t' \vect{v''} \in H,
\end{cases}$ $\Rarr$}

\noindent
$\Rarr$ $\{ \wtilde{p_j} \} = \frac{1}{\sin (t' - t'')} (\sin t' \vect{v''} - \sin t'' \vect{v'}) \in H$
($\sin (t' - t'') \ne 0$, because $t', t'' \in \mathbf{\Delta}$, $|\mathbf{\Delta}| \lesq \pi$), --- a contradiction.
So, $t'$ is the unique value providing the convergence of series in the left side of \equref{sphmainquadeq}.

We obtain $t'$ using the method analogous to that of Subcase 1a from Section 1.

$f_n(t) = \sum\limits_{j=1}^n (\wtilde{p_j} \cos t + \wtilde{q_j} \sin t)^2 - 1$ is non-decreasing relative to $n$, $f_n(t) \lesq f_{n+1}(t)$.
Rearranging, $f_n(t) = \alpha_n \cos^2 t + \beta_n \sin^2 t + \gamma_n \cos t \sin t - 1 = (\alpha_n - 1) \cos^2 t + (\beta_n - 1) \sin^2 t + \gamma_n \cos t \sin t$.

$\cos t \ne 0$, otherwise $\sum\limits_{j=1}^{\infty} \wtilde{q_j}^2 = 1$. Consequently, to solve $f_n(t) = 0$, we can divide it by $\cos t$:

\clin{$(\beta_n - 1) \tan^2 t + \gamma_n \tan t + (\alpha_n - 1) = 0$}

The quadratic trinomial here has no more than 2 roots ($\beta_n \rarr +\infty$ as $n\rarr \infty$), hence $f_n(t)$ has a finite number of zeroes in $\mathbf{\Delta}$.
Between them, the sign of the continuous $f_n(t)$ is constant. We consider it ``easy enough'' to determine, for each $n\in \mathbb{N}$,
the set $U_n = \{ t\in \mathbf{\Delta}\mid f_n(t) \lesq 0 \}$. It is clear that $U_n \supseteq U_{n+1}$.

Moreover, $\{ t' \} = \bigcap\limits_{n\in \mathbb{N}} U_n$; indeed, $f_n(t') \lesq f_{\infty}(t') = 0$,
and $\forall t \ne t'$ $\exists n_0$: $f_{n_0}(t) > 0$.

\textbf{Case 3:} $\sum\limits_{j=1}^{\infty} \wtilde{p_j}^2$ and $\sum\limits_{j=1}^{\infty} \wtilde{q_j}^2$ converge.
Then we denote $\alpha = \sum\limits_{j=1}^{\infty} \wtilde{p_j}^2$, $\beta = \sum\limits_{j=1}^{\infty} \wtilde{q_j}^2$,
$\gamma = 2 \sum\limits_{j=1}^{\infty} \wtilde{p_j} \wtilde{q_j}$, and rewrite \equref{sphmainquadeq} as
$(\alpha - 1) \cos^2 t + (\beta - 1) \sin^2 t + \gamma \sin t \cos t = 0$ $\LRarr$

\clin{$\LRarr$ $\left[\begin{array}{l}
\beta = 1, \quad \cos t = 0,\\
(\beta - 1) \tan^2 t + \gamma \tan t + (\alpha - 1) = 0, \quad \cos t \ne 0
\end{array}\right.$}

\textbf{Subcase 3a:} $(\beta - 1)^2 + \gamma^2 + (\alpha - 1)^2 > 0$. Then \equref{sphmainquadeq} again has a finite number of roots in $\mathbf{\Delta}$,
and they are ``easy to obtain'' (in accordance with our allowances). 

\textbf{Subcase 3b:} $\beta - 1 = \gamma = \alpha - 1 = 0$, therefore $\{ \wtilde{p_j} \} \in S$, $\{ \wtilde{q_j} \} \in S$, and $\{ \wtilde{p_j} \} \perp \{ \wtilde{q_j} \}$.
Then any $t\in \mathbf{\Delta}$ satisfies \equref{sphmainquadeq}.

Each root $t$, in turn, determines $\vect{s} = \{ \wtilde{p_j} \cos t + \wtilde{q_j} \sin t \}_j$.

This concludes the description of the solving method for SRP on $S$ with $d$.

{\ }

Again, there's the question about additional conditions for the set of sensors, or its extensions, providing the uniqueness of the solution,
especially in Subcase 3b with the most ``ambiguity''. We restrict our attention to the finiteness of the set of solutions.

{\ }

\prop{propExclSubcaseInf}
If $\vect{r}^{(2)} \perp \vect{r}^{(1)}$ and $\vect{r}^{(3)} \perp \vect{r}^{(1)}, \vect{r}^{(3)} \perp \vect{r}^{(2)}$, then Subcase 3b is impossible.

\begin{mproof}
In the basis $B$, not only $\vect{r}^{(1)} = (1; 0; 0; ...)$, but $\vect{r}^{(2)} = (0;1;0;0; ...)$ and $\vect{r}^{(3)} = (0;0;1;0;0; ...)$ then.
From \equref{sphcoeffcoordfromt} it follows that $\wtilde{p_k} = \cos t_k$ and $\wtilde{q_k} = \sin t_k$ for $k=1,2,3$;
therefore $\sum\limits_{j=1}^{\infty} (\wtilde{p_j}^2 + \wtilde{q_j}^2) \greq \sum\limits_{j=1}^3 (\wtilde{p_j}^2 + \wtilde{q_j}^2) = 3$.

Meanwhile, in Subcase 3b we have $\sum\limits_{j=1}^{\infty} \wtilde{p_j}^2 + \sum\limits_{j=1}^{\infty} \wtilde{q_j}^2 = \alpha + \beta = 2 < 3$.
\end{mproof}

The constraints for $\vect{r}^{(1)}$, $\vect{r}^{(2)}$, $\vect{r}^{(3)}$ here may be weakened:
$\sum\limits_{j=1}^3 (\wtilde{p_j}^2 + \wtilde{q_j}^2) > 2$ would suffice.

{\ }

\examp{exampDownDimSph} (analogous to \exampref{exampDownDimInfinite}).
\begin{mexamp}
Let $R$ be an orthonormal basis of $H$, $r^{(i)}_j = \delta_{ij}$, and let $\vect{s'} = \vect{r}^{(1)}$, $t' = 0$.
Then $t_1 = t' = 0$, $t_k = t' + d(\vect{r}^{(k)};\vect{s'}) = \arccos 0 = \frac{\pi}{2}$ for any $k \greq 2$.
$(\vect{s'};t')$ is the solution, on $S$, of the SRP defined by $(R;\{ t_i \}_{i\in \mathbb{N}})$.

Meanwhile, for any infinite $\widehat{R} \subset R$ such that $\vect{r}^{(1)} \notin \widehat{R}$, the truncated SRP
``$t_i = t + d(\vect{r}^{(i)};\vect{s})$ for $\vect{r}^{(i)} \in \widehat{R}$'' has no solution on $\widehat{S} = \{ \vect{x} \in L(\widehat{R}) \mid \| \vect{x} \| = 1 \}$.
\begin{mproof}
$\widehat{R}$ is the orthonormal basis of $L(\widehat{R})$, and we enumerate the elements of $\widehat{R}$ as they follow in $R$:
$\widehat{R} = \{ \hat{\vect{r}}^{(1)}; \hat{\vect{r}}^{(2)}; ...\}$; then, decomposing $\widehat{R}$ in itself, $\hat{r}^{(i)}_j = \delta_{ij}$.
Assuming $(\vect{s};t)$, $\vect{s} \in \widehat{S}$, to be the solution of truncated SRP, we have
$\frac{\pi}{2} \equiv \hat{t}_k = t + \arccos {<}\hat{\vect{r}}^{(k)};\vect{s}{>} = t + \arccos s_k$,
thus $s_k \equiv \cos(\frac{\pi}{2} - t) = \sin t$. Since $\dim L(\widehat{R}) = \infty$, $\vect{s} \in H$ implies $s_k \equiv 0$,
so $\vect{s} = \theta \notin \widehat{S}$, --- a contradiction.
\end{mproof}
\end{mexamp}

{\ }

{\ }

\textbf{Acknowledgements.} We would like to thank everyone who will provide the references to where these problems in $H$, ---
or perhaps, their more general forms, --- have been considered, studied, and solved already.

\end{document}